%Format: plain
\input amstex
\documentstyle{amsppt}

\input label.def
\input degt.def
\input dd.def
%\input debug.def
%\PrintLabels

\def\Dg:{\endgraf{\bf Dg:}\enspace\ignorespaces}

{\catcode`\@=11
\gdef\proclaimfont@{\sl}}

\loadbold

\def\todo:{{\bf to do:}}

\let\=\B
\def\ie{\emph{i.e\.}}
\def\eg{\emph{e.g\.}}
\def\cf.{\emph{cf\.}}
\def\all.{\emph{et al\.}}
\def\etc.{\emph{etc\.}}
\def\via{\emph{via}}

\newtheorem{Remark}\Remark\thm\endAmSdef

\def\dash{\item"\hfill--\hfill"}
\def\Dashes{\widestnumber\item{--}\roster}
\def\endDashes{\endroster}

\def\Cp#1{\Bbb P^{#1}}
\def\Rp#1{\Cp{#1}_\R}

\def\PGL{\operatorname{\text{\sl PGL}}}
\def\Qu{\operatorname{Qu}}
\def\Gr{\operatorname{Gr}}
\def\Sing{\operatorname{Sing}}
\def\Aut{\operatorname{Aut}}
\def\Pic{\operatorname{Pic}}
\def\ind{\operatorname{ind}}
\def\spcheck{\vee}

\let\Ga\alpha
\let\Gs\sigma
\let\Go\omega

\let\into\hookrightarrow

\def\tS{\tilde S}
\def\bX{\tilde X}

\def\CH{\Cal H}
\def\CM{\Cal M}

\def\go{\frak o}
\def\sd{_{\det}}

\def\bA{\bold A}
\def\bD{\bold D}
\def\bE{\bold E}
\def\bU{\bold U}

\def\SCH{S_{\CH}}

\def\fS{\frak S}
\def\fP{\frak P}

\def\PP{\Bbb P}

\def\sin{\text{\rm sp}_{\text{in}}}
\def\sout{\text{\rm sp}_{\text{out}}}

\def\tX{\tilde X}

\topmatter

\title
On real determinantal quartics
\endtitle

\author
Alex Degtyarev and Ilia Itenberg
\endauthor

\thanks
The second author is partially funded by
the ANR-09-BLAN-0039-01 grant
of {\it Agence Nationale de la Recherche}
and is a member of
%FRG "{\it Mirror Symmetry} \& {\it Tropical Geometry}".
FRG: Collaborative Research: Mirror Symmetry \& Tropical Geometry
(Award No. $0854989$).
\endthanks

\address
Department of Mathematics,\endgraf\nobreak
Bilkent University,\endgraf\nobreak
06800 Ankara, Turkey
\endaddress

\email
degt\@fen.bilkent.edu.tr
\endemail

\address
Universit\'e de Strasbourg,
IRMA and Institut Universitaire de France,
\endgraf\nobreak
7 rue Ren\'{e} Descartes\endgraf\nobreak
67084 Strasbourg Cedex, France
\endaddress

\email
ilia.itenberg\@math.unistra.fr
\endemail

\abstract
We describe all possible arrangements of the
ten nodes
of a generic real
determinantal quartic surface in $\Cp3$ with nonempty
spectrahedral region.
\endabstract



\endtopmatter

\document

\section{Introduction\label{S.intro}}

\subsection{Motivation}
It is a common understanding that, thank to the global Torelli
theorem for $K3$-surfaces~\cite{PSh-Sh} and surjectivity of the
period map~\cite{Kulikov}, any reasonable question concerning
the topology of
singular or real $K3$-surfaces can be reduced to a certain
arithmetic problem; many examples, treating the two subjects
separately, are found in the literature. However, there are but a
few papers where objects that are both real and singular are
considered; one can mention~\cite{It} and~\cite{Student}, which
deal, respectively, with real sextics with a single node in $\Cp2$
and real quartic surfaces with a single node in $\Cp3$.

In the present paper, we make an attempt to advance
this line of research, considering real quartic surfaces
with several nodes. Special attention is paid to
degenerations of nonsingular quartics, which are used
to control the topology of the resulting singular surfaces.
Since the classical problem of enumerating
all
equivariant
equisingular deformation types seems rather hopeless
(one would expect
thousands of classes),
we confine ourselves to a very special example
arising from convex algebraic geometry. Namely,
we describe arrangements of the
ten
nodes
of a generic
determinantal quartic with nonempty
spectrahedral region, see next subsection for details.

\subsection{Principal results}
Consider a generic dimension~$3$
real linear system~$V$ of quadrics in~$\Cp3$. Singular
quadrics form a surface $X\subset V\cong\Cp3$
of degree~$4$, which is
called a
\emph{transversal determinantal quartic}
(see Section~\ref{S.C}
for details and precise definitions).
In other words, we consider a
quartic surface $X\subset\Cp3$ given
by an equation of the form $\det\sum_{i=0}^3x_i\bar q_i=0$,
where $[x_0:x_1:x_2:x_3]$ are homogeneous coordinates in~$\Cp3$
and $\bar q_0,\bar q_1,\bar q_2,\bar q_3$ are
certain fixed nonzero symmetric
$(4\times4)$-matrices. Generically, such a surface is known to
have ten nodes.

Whenever present, quadrics given by definite quadratic forms
constitute a
single connected component of the complement $V_\R\sminus X_\R$;
this component is called the \emph{spectrahedral region} of~$V$.
This construction is a special case of a more general framework,
see~\cite{SS}, \cite{RS}
for details,
and the study of the shapes of various spectrahedra
is a major problem of convex algebraic geometry.

A transversal determinantal quartic has ten nodes, and the
original question posed to us by
B.~Sturmfels was whether all ten can be
located in the boundary of the spectrahedral region.
(The best known example, constructed explicitly,
had eight nodes in the boundary.)
We answer
this question in the affirmative; moreover, we describe all
possible arrangements of the nodes
with respect to the components of the complement
$V_\R\sminus X_\R$.

\theorem\label{th.main}
Let $X\subset\Cp3$ be a transversal real determinantal quartic
with nonempty spectrahedral region~$R$. Then $X$ has an
even number $m\ge0$ of real nodes in the boundary of~$R$ and an
even number $n\ge0$ of real nodes disjoint from~$R$, so
that $2\le m+n\le10$. Any pair of even numbers
$m,n\ge0$, $2\le m+n\le10$, is realized by a quartic as above.
\endtheorem

This theorem is proved in Subsection~\ref{proof.main}.

\Remark
It is worth emphasizing that any transversal
real determinantal quartic
with nonempty spectrahedral region has at least two real
nodes. Note that a similar, and even stronger, statement holds for
transversal real determinantal cubics in~$\Cp3$, which are
discriminants of linear systems of plane conics: such a cubic
(which is necessarily a Cayley cubic) has nonempty spectrahedral
region if and only if
%it has
%at least one real singular point.
at least one of its four nodes is real.
\endRemark

\subsection{Contents of the paper}
To prove Theorem~\ref{th.main}, we analyze the equisingular
stratification of the space of complex quartics in~$\Cp3$,
Section~\ref{S.quartics}, and identify the stratum
that is formed, up
to codimension one subset, by the transversal determinantal
quartics, Section~\ref{S.C}. Then we describe the sets of cycles
that can vanish under certain special nodal degenerations of a
real quartic surface, see Section~\ref{R.Q}. Finally, in
Section~\ref{S.R}, we show that each transversal real
determinantal quartic
is obtained by
a degeneration of a nonsingular quartic
with two nested spheres (the
so called \emph{hyperbolic quartic}),
and use previously known arithmetical
computations in order
to construct/prohibit various degenerations of the latter.

\subsection{Acknowledgements}
We are grateful to B.~Sturmfels for attracting our attention
to the problem and for motivating discussions.
This paper was conceived
during the second author's
stay at
Mathematical Sciences Research Institute at Berkeley
and completed during the first
author's visit to \emph{Universit\'e de Strasbourg}.

\section{Singular quartics in $\Cp3$\label{S.quartics}}

The principal result of this section is Theorem~\ref{enumeration},
which enumerates the equisingular strata of the space of
quartics in~$\Cp3$.

\subsection{Integral lattices}\label{s.lattice}
A {\it lattice\/} is a finitely generated free abelian group~$L$
supplied with a symmetric bilinear form $b:L\otimes L\to\ZZ$. We
abbreviate $b(x,y)=x \cdot y$ and $b(x,x)=x^2$. A lattice~$L$ is
{\it even\/} if $x^2=0\bmod2$ for all $x\in L$.
%The {\it
%determinant\/} $\det L$ is well defined.
As the transition matrix between two integral bases
has determinant $\pm1$, the
determinant $\det L\in\ZZ$
({\it i.e.}, the determinant
of the Gram matrix of~$b$ in
any
basis of~$L$)
is well defined.
A lattice~$L$ is called
{\it nondegenerate\/} if the determinant $\det L\ne0$; it is called {\it
unimodular\/} if $\det L=\pm1$.

Given a lattice~$L$,
the bilinear form can be extended to $L\otimes\Q$ by linearity. If
$L$ is nondegenerate, the dual group $L^\spcheck=\Hom(L,\ZZ)$ can
be identified with the subgroup
$$
\bigl\{x\in L\otimes\Q\bigm|
 \text{$x\cdot y\in\ZZ$ for all $x\in L$}\bigr\}.
$$
In particular,
%$L\subset L^\spcheck$.
$L$ itself can be identified with a subgroup
of $L^\spcheck$.

The group of isometries of a lattice~$L$ is denoted by
$\Aut L$.
Given a vector $a\in L$, $a^2\ne0$, the \emph{reflection} against
(the hyperplane orthogonal to)~$a$ is the automorphism
$r_a\:L\to L$, $x\mapsto x-2(x\cdot a)a/a^2$, provided that it is
well defined, \ie, takes integral vectors to integral vectors. The
reflection is always well defined if $a^2=\pm1$ or $\pm2$; if
$a^2=\pm4$, the reflection~$r_a$ is well defined if and only if
$a=0\bmod2L^\spcheck$.

A nondegenerate lattice~$L$ is called \emph{elliptic}
or \emph{hyperbolic} if its positive inertia index equals~$0$
or~$1$, respectively.
To any hyperbolic lattice~$H$ one can associate a hyperbolic space
$\PP(C):=C/\R^*$, where
$C=C_H:=\{x\in H\otimes\R\,|\,x^2>0\}$ is the
\emph{positive cone} of~$H$.
In particular, given a lattice $L$ and an isometric
involution $c\:L\to L$ with hyperbolic
invariant sublattice $L^c_+ = \{x\in L\,|\,c(x) = x\}$,
one can define the space $\PP(C^c_+)$.
Any subgroup $G\subset\Aut H$
generated by (some) reflections
%on~$H$
$r_a\:H\to H$
defined by vectors $a\in H$ with $a^2<0$
admits a polyhedral fundamental domain
$\fP_G\subset\PP(C)$:
it is the closure of (any) connected component
of the space $\PP(C)$ with all mirrors of~$G$ removed.

All lattices considered in the paper are even.
A \emph{root} in an
even lattice is a vector of square~$(-2)$.
A \emph{root system} is
an elliptic lattice generated by roots.
We use the standard notation $\bA_p$, $p\ge1$, $\bD_q$, $q\ge4$,
$\bE_6$, $\bE_7$, $\bE_8$ for the irreducible root systems of the
same name.
%We denote by $\bE_8$
%and~$\bA_1$ the negative definite lattices generated by the root
%systems of the same name,
%and we let
Let $\bU=\Z u_1\oplus\Z u_2$,
$u_1^2=u_2^2=0$, $u_1\cdot u_2=1$;
%the latter
this lattice is called
the \emph{hyperbolic plane},
and any basis $(u_1,u_2)$ as above is
called a \emph{standard basis} for~$\bU$.
Given a lattice~$L$ and an
integer~$d$, the notation $L(d)$
stands for the lattice obtained from~$L$ by multiplying the values
of the bilinear form by~$d$.

%\subsection{The singular homological type}
\subsection{Singular homological types}

\definition\label{def.Sigma}
A \emph{set of \rom(simple\rom) singularities} is a pair
$(\Sigma,\Gs)$, where $\Sigma$ is a root system and
$\Gs$ is a collection of roots of~$\Sigma$ constituting a Weyl
chamber of~$\Sigma$. An isometry $\Sigma_1\to\Sigma_2$ of two
sets of singularities
$(\Sigma_i,\Gs_i)$, $i=1,2$, is \emph{admissible} if it
takes~$\Gs_1$ to~$\Gs_2$.
\enddefinition

\Remark
Any Weyl chamber of a root system~$\Sigma$ can be taken to any
other Weyl chamber by an element of the Weyl group of~$\Sigma$,
which extends to any larger lattice containing~$\Sigma$. For this
reason, when speaking about the isomorphism classification of sets
of singularities, configurations, and singular homological types (see
below), the subset~$\Gs$ in Definition~\ref{def.Sigma} can be, and
often is disregarded.
\endRemark

\definition\label{configuration}
A \emph{configuration} (extending a given set of singularities
$(\Sigma,\Gs)$) is a finite index extension
$\tS\supset S:=\Sigma\oplus\Z h$, $h^2=4$,
satisfying the following conditions:
\roster
\item\local{conf.1}
each root $r\in\tS\cap(\Sigma\otimes\Q)$ belongs
to~$\Sigma$;
\item\local{conf.2}
$\tS$ does not contain an element~$u$ with $u^2=0$ and
$u\cdot h=2$.
\endroster
An \emph{admissible isometry} of two configurations
$\tS_i\supset S_i=\Sigma_i\oplus\Z h_i$, $i=1,2$, is an isometry
$\tS_1\to\tS_2$
taking~$h_1$ to~$h_2$ and inducing an admissible isometry
$\Sigma_1\to\Sigma_2$.
\enddefinition

\definition\label{def.AHT}
A \emph{singular homological type} (extending a set
of singularities $(\Sigma,\Gs)$) is an extension of the
orthogonal direct sum
$S:=\Sigma\oplus\Z h$, $h^2=4$, to a lattice~$L$ isomorphic to
$2\bE_8\oplus3\bU$, such that the primitive hull~$\tS$ of~$S$ in~$L$
is a configuration.
(The singular homological type is also said to \emph{extend} the
configuration $\tS\supset S$.)
An \emph{isomorphism} between two singular homological types
$L_i\supset S_i\supset\Gs_i\cup\{h_i\}$, $i=1,2$, is an
isometry
$L_1\to L_2$ taking~$h_1$ to~$h_2$ and $\Gs_1$ to~$\Gs_2$ (as a set).
\enddefinition

A singular homological type is uniquely determined by
the collection $\CH=(L,h,\Gs)$; then
$\Sigma=\Sigma_{\CH}$ is the sublattice spanned
by~$\Gs$, and $S=\SCH=\Sigma\oplus\Z h$.

Given a singular homological type~$\CH$, the orthogonal complement
$\SCH^\perp$ is a nondegenerate lattice of positive inertia
index~$2$.
Hence, the orthogonal projection of any
positive definite $2$-subspace $\Go_1\subset\SCH^\perp\otimes\R$
to any other such subspace $\Go_2$ is an isomorphism of vector
spaces; it can be used to compare orientations of~$\Go_1$
and~$\Go_2$. Thus, a choice of an orientation of one positive
definite $2$-subspace in $\SCH^\perp\otimes\R$ defines a coherent
orientation of any other.

\definition\label{def.orientation}
An \emph{orientation} of a singular homological type
$\CH=(L,h,\Gs)$ is a choice of coherent orientations of
positive definite $2$-subspaces of $\SCH^\perp\otimes\R$.
Oriented singular homological types $(\CH_i,\go_i)$, $i=1,2$, are
\emph{isomorphic} if there is an isomorphism $\CH_1\to\CH_2$
taking~$\go_1$ to~$\go_2$.
A singular homological type~$\CH$ is called
\emph{symmetric} if $(\CH,\go)\cong(\CH,-\go)$, \ie, it
$\CH$ admits an automorphism reversing
orientation.
\enddefinition

\subsection{Classification of singular quartics}
Let $X\subset\Cp3$ be a quartic surface with simple singularities
only. Denote by $\bX\to X$ the minimal resolution of singularities
of~$X$; it is a minimal $K3$-surface. Introduce the following
objects:
\Dashes

\dash
$L_X=H_2(\bX)=H^2(\bX)$,
regarded as a lattice \via\ the intersection form
(we always identify homology and cohomology \via\ the Poincar\'e
duality);

\dash
$\Gs_X\subset L_X$, the set of the classes of the
exceptional divisors contracted by the blow-up map $\bX\to X$;

\dash
$h_X\in L_X$, the class of the pull-back of a generic plane section
of~$X$;

\dash
$\Go_X\subset L_X\otimes\R$, the oriented $2$-subspace
spanned by the real and imaginary parts of the class of a
holomorphic $2$-form on~$\bX$ (the \emph{period} of~$\bX$).

\endDashes
Note that $\Go_X$ is positive definite.
According to~\cite{Saint-Donat;Urabe}, a triple $\CH=(L,h,\Gs)$ has the form
$(L_X,h_X,\Gs_X)$ for a quartic $X\in\Cp3$ as above if and only if it
is a singular homological type in the sense of Definition~\ref{def.AHT}.
%In this case,
If this is the case,
the above orientation of $\Go_X$
defines an orientation of~$\CH$.

The following
theorem is
quite expectable;
however, we could not find an explicit statement in the
literature. The surjectivity part is contained in~\cite{Urabe}.

\theorem\label{enumeration}
The map sending a quartic surface $X\subset\Cp3$
with simple singularities to the pair
consisting of its singular homological type $\CH_X=(L_X,h_X,\Gs_X)$ and the
orientation of the space~$\Go_X$ establishes a one-to-one
correspondence between the set of equisingular deformation classes of
quartics
with a given set of simple singularities $(\Sigma,\Gs)$ and the set of
isomorphism classes of oriented abstract singular homological types
extending~$(\Sigma,\Gs)$.
\endtheorem

\proof
Proof of this theorem repeats, almost literally, the proof of a
similar theorem for plane sextic curves, see~\cite{JAG}. It is
based on Beauville's construction~\cite{Beauville}
of a fine period space of marked polarized $K3$-surfaces. We omit
the details.
\endproof

The equisingular stratum of the space of quartic surfaces
in~$\Cp3$ corresponding to
an oriented singular homological type $(\CH,\go)$ will be denoted by
$\CM(\CH,\go)$. If $\CH$ is symmetric, we abbreviate this notation
to $\CM(\CH)$.
As part of the proof of Theorem~\ref{enumeration}, one obtains an
explicit description of the moduli space of quartics, which
results in the following formula for its dimension
$$
\dim\CM(\CH,\go)/\!\PGL(4,\C)=19-\rank\Sigma_{\CH}
\eqtag\label{eq.dim}
$$
(similar to the corresponding formula for plane sextics).
Note that $\rank\Sigma_{\CH}=\#\Gs$ equals the total Milnor
number~$\mu(X)$ of~$X$.

\Remark
The equisingular deformation classification of quartic surfaces
with isolated singularities and at least one non-simple singular
point is found in~\cite{quartics;quartics.2}. With a few
exceptions, the deformation class of such a quartic is also
determined by its (appropriately defined) singular homological type.
\endRemark

\section{Real quartics\label{R.Q}}

In this section, we analyze the position of the vanishing cycles
of a degeneration of a nonsingular quartic with respect to its
period domain. The principal results are Theorems~\ref{th.only.if}
and~\ref{th.if}.

\subsection{Real homological types}
Given an isometric
involution $c\:L\to L$ on a lattice~$L$, we denote
by
$L^c_\pm=\{x\in L\,|\,c(x)=\pm x\}\subset L$
the $(\pm1)$-eigenlattices of~$c$.
If $L$ is nondegenerate, $L^c_\pm$ are the orthogonal complements
of each other.

\definition\label{def.NRHT}
A
\emph{real homological type}
is a triple $(L, h, c)$, where
$L$ is a lattice isomorphic to~$2\bE_8\oplus 3\bU$,
$h \in L$
is a vector of square~$4$,
and $c\: L \to L$ is an isometric involution
such that
\Dashes
\dash
the sublattice $L^c_+$ is hyperbolic, and
\dash
one has $h \in L^c_-$.
\endDashes
An \emph{isomorphism} between two
real homological types
$(L_i, h_i, c_i)$, $i=1,2$, is an
isometry
$\varphi\: L_1\to L_2$ such that $\varphi(h_1) = h_2$
and $\varphi \circ c_1 = c_2 \circ \varphi$.
\enddefinition

\definition\label{def.poly}
Let $(L, h, c)$ be a
real homological type.
Consider vectors $e \in L^c_+$ of the following three kinds:
\roster
\item\local{s2}
$e^2=-2$, \ie, $e$ is a root;
\item\local{s4h}
$e^2=-4$ and $e=h\bmod2L$;
\item\local{s4}
$e^2=-4$, $e\ne h\bmod2L$, and~$e$ is
\emph{decomposable},
\ie, $e=r'-r''$ for a pair of
%orthogonal
roots $r',r''\in L$ such that
$r'\cdot r''=r'\cdot h=r''\cdot h=0$ and
$r''=-c(r')$.
\endroster
A \emph{fundamental tower} of $(L, h, c)$
is a
triple $\fS\subset\fP\subset\bar\fP\subset\PP(C^c_+)$ of
%polyhedral
fundamental
domains
of the subgroups of $\Aut L^c_+$
generated by the reflections defined by all vectors of~$L^c_+$ of
type \loccit{s2}--\loccit{s4},
\loccit{s2}--\loccit{s4h}, and~\loccit{s2}, respectively.
\enddefinition

Note
that, in cases~\loccit{s4h} and~\loccit{s4}, the conditions
imposed imply
%that
$e=0\bmod2(L^c_+)^\spcheck$, \ie, $e$ does
define a reflection $r_e\:L^c_+\to L^c_+$.
Moreover, one can easily see that
this reflection extends to an automorphism of the homological
type. As a consequence,
any two fundamental towers are related by an automorphism of the
homological type (in fact, by a sequence of reflections).

\definition\label{def.periods}
A
%nonsingular
real homological type $(L,h,c)$ equipped with a
distinguished fundamental tower
$\fS\subset\fP\subset\bar\fP$ is called a
\emph{period lattice}, and the polyhedra~$\fP$
and~$\bar\fP$ are called the \emph{period domains} (more precisely,
the period domain of real quartics and
that of abstract real $K3$-surfaces, respectively).
\enddefinition

%Any two fundamental towers are related by an automorphism of the
%homological type (in fact, by a sequence of reflections).
As explained in Subsection~\ref{s.lattice},
the facets of the polyhedra
in Definition~\ref{def.poly} are (parts of)
some
of the mirrors (walls)
of the respective groups, \ie, hyperplanes orthogonal
to vectors of corresponding types. We refer to the type of the
vector as the \emph{type} of the corresponding wall.

\Remark\label{rem.fS}
The polyhedra~$\fP$ and $\bar\fP$ have a
certain geometric meaning (see Subsection~\ref{real-periods} below),
whereas $\fS$ does not.
%Furthermore, a fundamental tower is unique
%up to automorphism of the homological type, whereas a refinement
%is \emph{a priori} not.
However, in
many examples, $\fS$ is much easier to compute
and, on the other hand, a choice of~$\fS$
determines the other two polyhedra: $\fP$ is paved by the copies
of~$\fS$ obtained from~$\fS$ by iterated reflections against
(the consecutive images of) the walls of type~\iref{def.poly}{s4},
and, similarly, $\bar\fP$ is paved by the copies of~$\fP$ obtained
by iterated reflections against the walls of
type~\iref{def.poly}{s4h}.
\endRemark

\subsection{Invariant periods}\label{real-periods}
A quartic $X\subset\Cp3$
is called \emph{real} if it is invariant under the complex
conjugation involution $\conj\:\Cp3\to\Cp3$. The
involution~$\conj$ restricts to~$X$ and, if $X$ is singular, lifts
to the minimal resolution~$\bX$ of singularities of~$X$, turning
both into \emph{real $K3$-surfaces}.
A nonsingular real quartic $X\subset\Cp3$ gives rise to a
real homological
type $(L_X, h_X, c_X)$, where $c_X\: L_X \to L_X$
is the involution induced by $\conj$.

Two nonsingular
real quartics $X,Y\subset\Cp3$ are said to be \emph{coarse
deformation equivalent} if $X$ is equivariantly deformation
equivalent to either~$Y$ or the quartic $Y'$ obtained from~$Y$ by
an orientation reversing automorphism of~$\Cp3$. A coarse
deformation class consists of one or two components of the space
of nonsingular real quartics; in the former case, the quartics are
called \emph{amphichiral}, in the latter case, \emph{chiral}.
The following statement is found in~\cite{Nikulin}.

\theorem\label{th-nikulin}
Two nonsingular real quartics $X,Y \subset \Cp3$ are coarse
deformation equivalent if and only if
the corresponding
real homological types
$(L_X, h_X, c_X)$ and $(L_Y, h_Y, c_Y)$ are isomorphic.
\qed
\endtheorem

A complete classification of nonsingular
real quartics in~$\Cp3$
up to equivariant deformation, addressing in
particular the chirality problem, and an
interpretation of the result in topological terms
are found
in~\cite{Kharlamov}.

Given a real quartic $X\subset\Cp3$ or, more generally, a real
$K3$-surface $(X,\conj)$,
a holomorphic $2$-form $\Omega_X$ on~$X$ can be normalized
(uniquely up to a nonzero real factor) so that
$\conj^*\Omega_X=\bar\Omega_X$; such a form
is called \emph{real}.
The real part $(\omega_X)_+$
of the class of a real form $\Omega_X$
belongs to $(L_X)^{c_X}_+ \otimes \R$
and defines
a point $[(\omega_X)_+]$
in the associated hyperbolic space;
this point
is called the \emph{invariant period} of~$X$.

%Let $X\subset\Cp3$ be a nonsingular real quartic, and
Fix a period lattice $(L,h,c;\fS\subset\fP\subset\bar\fP)$ with
the
%underlying
real homological type $(L,h,c)$
isomorphic to that of~$X$. A
particular choice of an isomorphism
$\varphi\:(L_X,h_X,c_X)\to(L,h,c)$ is called a \emph{marking}
of~$X$. A marking~$\varphi$ is called \emph{proper} if
%$\varphi[\omega_X]_+\in\fP$.
$\varphi[(\omega_X)_+]\in\fP$.
In fact,
if~$X$ is nonsingular, the image
%$\varphi[\omega_X]_+$
$\varphi[(\omega_X)_+]$
under a proper marking belongs to the
interior $\Int\fP$, see, \eg,~\cite{Nikulin}. It follows that any
two proper markings differ by a symmetry of~$\fP$.

%\subsection{Degenerations of real quartics}\label{degenerations}
\subsection{Degenerations}\label{degenerations}
A \emph{degeneration} is a smooth family $X_t\subset\Cp3$,
$t\in[0,1]$, of
real quartics such that all quartics $X_t$, $t\in(0,1]$ are
nonsingular. For simplicity, we confine ourselves to the case when
$X_0$ has simple nodes only as singularities. Recall that the
homology groups $H_2(X_t)$ of the nonsingular members of the
family are canonically identified \via\ the Gauss-Manin
connection, and this common group contains a set of
\emph{vanishing cycles} (defined up to sign), one for each node
of~$X_0$.

The Gauss-Manin connection can be extended to identify the
homology of~$X_1$ with the homology of the minimal
resolution~$\tilde X_0$ of~$X_0$, taking (up to sign)
the vanishing cycles to
the classes of the exceptional divisors contracted
in~$X_0$.

At each real node of~$X_0$, the difference of the local Euler
characteristics of the real parts of~$X_0$ and~$X_1$ is $\pm1$;
according to this difference, the node is called \emph{positive}
or \emph{negative}, respectively. Negative are the nodes whose
vanishing cycles are $c_{X_1}$-invariant. Below, we are interested
in the \emph{non-positive} nodal degenerations, \ie, such that
each node of $X_0$ is either not real or real and negative.

\definition\label{def.cycles}
Let $(L,h,c;\fS\subset\fP\subset\bar\fP)$ be a period lattice.
A collection of roots
$r_i,s'_j,s''_j\in L$, $i=1,\ldots,k$, $j=1,\ldots l$,
is called an \emph{admissible system of
cycles}
if it satisfies the following conditions:
\roster
\item\local{vc.perp}
all roots are orthogonal to each other and to~$h$;
\item\local{vc.no.root}
the primitive hull in~$L$ of the sublattice spanned by $r_i$,
$s'_j$, $s''_j$ contains no roots other than $\pm r_i$,
$\pm s'_j$, $\pm s''_j$, \cf.~\iref{configuration}{conf.1};
\item\local{vc.fP}
each root $r_i$, $i=1,\ldots,k$, belongs to $L^c_+$ and
defines a facet
of~$\fP$, which is necessarily of type~\iref{def.poly}{s2};
\item\local{vc.fS}
for each $j=1,\ldots,l$, one has $c(s'_j)=-s''_j$
and the decomposable
invariant vector $s'_j-s''_j$ defines a
type~\iref{def.poly}{s4} facet of~$\fS$.
\endroster
\par\removelastskip
\enddefinition

\theorem\label{th.only.if}
Let
$(L,h,c;\fS\subset\fP\subset\bar\fP)$ be a period lattice,
and let $X_t$ be a non-positive nodal degeneration with the real
homological type of the nonsingular
surface $X:=X_1$ isomorphic to $(L,h,c)$.
%$X \subset \Cp3$ a nonsingular real quartic
%whose real homological type is isomorphic
%to $(L, h, c)$,
%and $X_t$ a non-positive nodal degeneration
%with $X = X_1$.
Then
%the surface
$X$ admits a proper marking that takes
the set of vanishing cycles of~$X_t$
%the degeneration
to an
admissible system of cycles.
\endtheorem

\theorem\label{th.if}
Given a period lattice $(L,h,c;\fS\subset\fP\subset\bar\fP)$ and
an admissible system of cycles $\Gs=\{r_i,s'_j,s''_j\}$,
$i=1,\ldots,k$, $j=1,\ldots,l$, there exists a nonsingular real
quartic~$X$ and a proper marking
$\varphi\:(L_X,h_X,c_X)\to(L,h,c)$ which identifies~$\Gs$ with the set of
vanishing cycles of a certain non-positive nodal degeneration
of~$X$.
\endtheorem

\proof[Proof of Theorem~\ref{th.only.if}]
Clearly, the vanishing cycles are orthogonal to each other and
to~$h$, \ie, satisfy~\iref{def.cycles}{vc.perp},
as they are geometrically disjoint and can be chosen
disjoint from a hyperplane section.

Denote by $r_i\in L_X$, $i=1,\ldots,k$, the vanishing cycles
corresponding to the real nodes of~$X_0$, and by $s'_j,s''_j\in L_X$,
$j=1,\ldots,l$, those corresponding to the pairs of complex
conjugate nodes; the latter are oriented so that
$c_X(s'_j)=-s''_j$.

%Consider a degeneration~$X_t$
%with $X=X_1$.
Let
$\tX_0$ be the minimal
resolution of singularities of~$X_0$.
Recall that,
using the Gauss-Manin connection,
we identify the homology of~$X$ and~$\tX_0$.
%, let
Let
$\tilde\omega\in L_X\otimes\C$ be the class realized by
a real holomorphic $2$-form on~$\tX_0$,
and let $[\tilde\omega_+]\in\fP_X$
be its invariant part (the invariant period of~$\tX_0$).
Note that $[\tilde\omega_+]$ does belong to~$\fP_X$ as it is the
limit of invariant periods of~$X_t$, which are all
in $\Int\fP_X$.
One has
$\Pic\tX_0=\tilde\omega^\perp\cap L_X$. Denote
%$\Pic_h\tX_0=h^\perp$,
$\Pic_h\tX_0 = (h_X)^\perp$,
the orthogonal complement of~$h_X$
%$h$
in $\Pic\tX_0$.
Up to sign, any root in
$\Pic_h\tX_0$
is represented
by a unique $(-2)$-curve contracted in $X_0$, and these are all
$(-2)$-curves contracted. It follows that the roots
of
$\Pic_h\tX_0$
are precisely the vanishing cycles;
in particular, this implies
condition~\iref{def.cycles}{vc.no.root}.
Thus, the maximal root system in $\Pic_h\tX_0$
is $(k + 2l)\bA_1$,
and all its roots define a common face
of all its Weyl chambers. Passing to the $c_X$-invariant part,
one easily concludes that the invariant vanishing cycles~$r_i$,
$i=1,\ldots,k$,
define a
common face of all $\fP$-like fundamental polyhedra
containing
$[\tilde\omega_+]$,
in particular, of~$\fP_X$, whereas the
decomposable vectors $s_j'-s_j''$, $j=1,\ldots,l$, define a common
face of
all $\fS$-like polyhedra
containing
$[\tilde\omega_+]$;
for
the latter, one can take any polyhedron
$\fS'$
containing $[\tilde\omega_+]$ and contained
in~$\fP_X$.
Due to Theorem~\ref{enumeration} and
condition~\iref{configuration}{conf.2} in the definition, one has
$s'_j-s''_j\ne h\bmod2L$; hence the wall defined by this vector is
of type~\iref{def.poly}{s4}.

It remains to consider any proper marking of~$X_1$
and, if necessary, adjust it by a symmetry
of~$\fP$ to make sure that the polyhedron
%~$\fS_X$
$\fS'$
constructed above
is taken to the preselected polyhedron~$\fS$.
%This completes the
%proof of Theorem~\ref{th.only.if}.
\endproof

\proof[Proof of Theorem~\ref{th.if}]
Let $f_{\fP}$ and $f_{\fS}$ be the intersections of
the facets defined
in~\iref{def.cycles}{vc.fP} and~\ditto{vc.fS}, respectively.
Notice that $f_{\fP}$ and $f_{\fS}$ are nonempty faces
of $\fP$ and $\fS$, respectively
(since the facets intersected are mutually orthogonal).
One
has $f_{\fP}\perp f_{\fS}$ and, since the symmetry
about~$f_{\fS}$ preserves~$\fP$, it also preserves~$f_{\fP}$. It
follows that the subspace supporting~$f_{\fS}$ intersects
$f_{\fP}$ at at least one interior point. Let $[\tilde\omega_+]$
be such a point,
and let $[\tilde\omega_-]\subset\PP(C^c_{-h})$
be a
point in the intersection of the hyperplanes defined
by the skew-invariant vectors $s'_j+s''_j$, $j=1,\ldots,l$, in the
hyperbolic space associated with the orthogonal
complement~$L^c_{-h}$ of~$h$ in~$L^c_-$. (Since all hyperplanes
are orthogonal to each other, they obviously intersect.)
Due to~\iref{def.cycles}{vc.no.root},
the pair $([\tilde\omega_+], [\tilde\omega_-])$ can be chosen
generic in the sense that $[\tilde\omega_+]$
and $[\tilde\omega_-]$ are not simultaneously orthogonal to
any root of~$L$ which is orthogonal to~$h$ and does not belong to
$\pm\Gs$.

Let $U \subset \PP(C^c_{-h})$
be a sufficiently small neighborhood of $[\tilde\omega_-]$.
Consider a generic path
$([(\omega_+)_t],[(\omega_-)_t])\in\Int\fP \times U$,
$t\in(0,1]$,
converging to the point $([\tilde\omega_+],[\tilde\omega_-])$.
According to~\cite{Nikulin}, it gives rise to a family~$X_t$ of
properly marked nonsingular
real quartics. (Strictly speaking, the path used should avoid a certain
codimension~$2$ subset, see~\emph{loc\. cit\.} for the technical
details.) This family can be chosen
to converge to a singular quartic~$X_0$ (\cf.,
\eg,~\cite{Shimada} for a detailed proof for the similar case of
plane sextic, \ie, polarization of square~$2$),
%which
and the limit quartic~$X_0$
is necessarily
real.
As in the previous proof,
considering the Picard group $\Pic_h\tX_0$ and using the fact that the
pair $([\tilde\omega_+],[\tilde\omega_+])$ is generic, one concludes
that the irreducible $(-2)$-curves contracted in~$X_0$ are
precisely those realizing the elements of~$\Gs$ (here, crucial is
condition~\iref{def.cycles}{vc.no.root}, which implies that
$(L,h,\Gs)$ is a singular homological type); hence, these
elements are the vanishing cycles.
\endproof

\Remark
Note that,
if negative nodes are present, the real
structure does \emph{not} change continuously on the
desingularized family~$\tX_t$ of abstract $K3$-surfaces;
in fact, the real homological type of
the limit surface~$\tX_0$, defined in the
obvious way, is not even isomorphic to that of~$X_t$, $t>0$.
However, the real structure does
change continuously on the quartics.
\endRemark

\section{Complex determinantal quartics\label{S.C}}

The goal of this section is Theorem~\ref{Hdet},
which identifies the
equisingular stratum containing transversal determinantal
quartics.

\subsection{Notation}
Let $\Qu(n)\cong\Cp{N(n)}$ be the space of quadrics in
$\Cp{n}$; here $N(n)=\frac12n(n+3)$. Let, further,
$\Qu_r(n)\subset\Qu(n)$,
$0\le r\le n$, be the space of quadrics of corank~$r$. The closure
$\Delta(n)$ of $\Qu_1(n)$ is called the \emph{discriminant
hypersurface}; it has degree $n+1$.

The singular locus of a quadric~$Q$ of corank $r>0$ is an
$(r-1)$-subspace of~$\Cp{n}$. Sending~$Q$ to its singular locus,
one obtains a locally trivial fibration
$$
\Qu_r(n)\to\Gr(n+1,r);
\eqtag\label{eq.fibration}
$$
its fiber is $\Qu_0(n-r)$. (We let $\Qu(0)=\Qu_0(0)=\pt$.)
Thus, $\Qu_r(n)$ is a smooth
quasi-projective variety and
$\dim\Qu_r(n)=N(n)-\frac12r(r-1)$.

\definition
A \emph{geometric hyperplane} is a hyperplane
$$
H_p:=\bigl\{Q\in\Qu(n)\bigm|Q\ni p\bigr\}
$$
consisting of all
quadrics passing through a fixed point $p\in\Cp{n}$.
\enddefinition

For each point $p\in\Cp{n}$,
fibration~\eqref{eq.fibration} restricts to
a locally trivial fibration
$$
\Qu_r(n)\sminus H_p\to\Gr(n+1,r)\sminus\Gr(n,r-1)
\eqtag\label{eq.fibration.2}
$$
with fiber $\Qu_0(n-r)\sminus H_{p'}$. Here, the difference in the
right hand side is the space of all $(r-1)$-planes in~$\Cp{n}$ not
passing through~$p$.

Since, in this paper, we are mainly concerned with quadrics
in~$\Cp3$, we abbreviate the notation as follows: let
$\Qu=\Qu(3)$, $\Delta=\Delta(3)$, and let
$\Delta'$ and $\Delta''$ be the closures of $\Qu_2(3)$ and
$\Qu_3(3)$, respectively. Let also
$\Delta^\circ=\Qu_1(3)=\Delta\sminus\Delta'$. One has
$$
\dim\Qu=9,\quad
\dim\Delta=\dim\Delta^\circ=8,\quad
\dim\Delta'=6,\quad
\dim\Delta''=3.
$$
Recall also that $\deg\Delta=4$ and $\deg\Delta'=10$
(see~\cite{HT}).
%\mnote{It:
%reference added; this article contains a general formula,
%do we need it?\Dg: I think not. Is it long?
%\endgraf It: the answer is
%$\prod_{i = 0}^{r - 1}(C^{r - i}_{n + i})/(C^{i}_{2i + 1})$,
%where $r$ is a corank; I also do not think that we should
%include the formula}

Let $V$ be a subspace of~$\Qu$ of dimension~$3$.
Unless $V\in\Delta$, the
intersection $\Delta_V:=V\cap\Delta$ is a quartic
in~$V$. Any quartic $X\in\Cp3$ such that the pair
$(\Cp3,X)$ is isomorphic to $(V,\Delta_V)$ above is called a
\emph{determinantal quartic}.

A $3$-space~$V$ is called \emph{transversal} if it is
transversal to the strata $\Delta^\circ$,
$\Delta'\sminus\Delta''$, and~$\Delta''$. Any determinantal
quartic $X\subset\Cp3$ isomorphic to $\Delta_V\subset V$ is also
called \emph{transversal}.
If $V$ is transversal,
the singular locus $\Sing\Delta_V$ coincides with
$V\cap\Delta'$ and consists of ten type~$\bA_1$ points.
Conversely, if $\Sing\Delta_V$ consists of ten type~$\bA_1$
points, $V$ is transversal.

For a $3$-space $V\subset\Qu$ as above,
we denote $\Delta^\circ_V=\Delta_V\sminus\Delta'$.

\subsection{Some fundamental groups}
Observe that
the set $\Sing\Delta(n)=\Qu_{\ge2}(n)$ of singular points
of~$\Delta(n)$ has
codimension~$3$ in $\Qu(n)$. Hence
a generic plane section of $\Delta(n)$ is a
nonsingular plane curve of degree $n+1$ and, due to Zariski's
hyperplane section theorem~\cite{Zariski}, one has
$\pi_1(\Qu_0(n))=\Z_{n+1}$. A generic plane section of the union
$\Delta(n)\cup H_p$ is a transversal union of a nonsingular curve
and a line; hence, $\pi_1(\Qu_0(n)\sminus H_p)=\Z$.

\proposition\label{pi1}
One has $\pi_1(\Delta^\circ)=0$.
\endproposition

\proof
The Serre exact sequence of
fibration~\eqref{eq.fibration} takes the form
$$
\CD
\pi_2(\Cp3)@>>>\pi_1(\Qu_0(2))@>>>\pi_1(\Delta^\circ)@>>>\pi_1(\Cp3)\\
@|@|@.@|\\
\Z@.\Z_3@.@.0.
\endCD
$$
From this sequence, one concludes that
$\pi_1(\Delta^\circ)=H_1(\Delta^\circ)$ is a
quotient of~$\Z_3$ and, moreover, the inclusion homomorphism
$H_1(U\sminus\Delta')\to H_1(\Delta^\circ)$ is onto, where $U$ is
a regular neighborhood in~$\Delta$ of a point $Q\in\Delta'$. A
normal $3$-plane
section of~$\Delta'$ in~$\Delta$ is a type~$\bA_1$
singularity and $U\sminus\Delta'$ is homotopy equivalent to its
link.
Hence, one has $H_1(U\sminus\Delta')=\Z_2$ and the statement
follows.
\endproof

\proposition\label{pi1.affine}
One has $\pi_1(\Delta^\circ\sminus H_p)=\Z_2$
\rom(for any point $p\in\Cp3$\rom).
\endproposition

\proof
Similar to the previous proof, using
fibration~\eqref{eq.fibration.2} instead of~\eqref{eq.fibration},
one concludes that the abelian group
$\pi_1(\Delta^\circ\sminus H_p)=H_1(\Delta^\circ\sminus H_p)$
is a quotient of the group
$H_1(U\sminus\Delta')=\Z_2$, where $U$ is
a regular neighborhood in~$\Delta$ of a point in~$\Delta'$. Thus,
it remains to show that $\Delta^\circ\sminus H_p$ admits a
nontrivial double covering.

Let $Q\in\Delta^\circ$, $Q\not\ni p$. Denote by~$q$ the (only)
singular point of~$Q$. Then, there are exactly two planes passing
through the line $(pq)$ and tangent to~$Q$ along a whole
generatrix: the original quadric~$Q$ and the two planes are the
cones, with the vertex at~$q$, over the section of~$Q$ by a
generic plane $\Ga\ni p$ and the two tangents
%in~$P$
to this section
passing through~$p$.
Clearly, the space of all pairs
$$
(Q,\{\text{tangent plane as above}\})
$$
is a double covering of $\Delta^\circ\sminus H_p$. This covering
is nontrivial: for example, in the family
$$
Q_t=\bigl\{(x_1-x_3)^2+x_1^2-e^{2\pi it}x_2^2=0\bigr\},
\qquad t\in[0,1],
$$
the two tangents $x_1=\pm e^{\pi it}x_2$ are interchanged. (In
particular, it follows that the path $Q_t$, $t\in[0,1]$, is a
non-contractible loop in
$\Delta^\circ\sminus H_{(0:0:0:1)}$.)
\endproof

\Remark
Similar to Propositions~\ref{pi1} and~\ref{pi1.affine}, one can
easily show that all fundamental groups $\pi_1(\Qu_r(n))$ and
$\pi_1(\Qu_r(n)\sminus H_p)$ are cyclic.
\endRemark

\corollary\label{cor.pi1}
Let $p\in\Cp3$. Then, for a
generic transversal $3$-plane $V\subset\Qu$, one
has
$\pi_1(\Delta^\circ_V)=0$ and
$\pi_1(\Delta^\circ_V\sminus H_p)=\Z_2$.
\endcorollary

\proof
The statement follows from Propositions~\ref{pi1}
and~\ref{pi1.affine} and Zariski type hyperplane section theorem
for quasi-projective varieties (see~\cite{Hamm;GM1;GM2} or recent
survey \cite{Cheniot, Theorem~5.1}).
\endproof

\subsection{The determinantal stratum}
In this subsection, we identify the stratum in the space of
quartics formed by the transversal determinantal ones.

\lemma\label{non-coplanar}
Any transversal determinantal quartic has a quadruple of
non-coplanar singular points.
\endlemma

We postpone the proof of this technical statement till
Subsection~\ref{proof.non-coplanar}.

\lemma\label{irreducible}
The space of quintuples $(V;Q_1,Q_2,Q_3,Q_4)$, where $V\subset\Qu$
is a transversal $3$-space and $Q_1$, $Q_2$, $Q_3$, $Q_4$ are four
non-coplanar singular points of $\Delta_V$, is an irreducible
quasi-projective variety of dimension~$24$.
\endlemma

\proof
The statement is a tautology, as the $3$-subspace
$V\subset\Qu$ is uniquely determined by
a quadruple of its non-coplanar
points $Q_1,Q_2,Q_3,Q_4$.
Thus, the space in question is a Zariski open
subset of the irreducible variety $(\Delta'\sminus\Delta'')^4$.
\endproof

Comparing the dimensions, see~\eqref{eq.dim},
one arrives at the following corollary.

\corollary\label{whole.stratum}
Transversal determinantal quartics $X\subset\Cp3$ form a Zariski
open
subset of a single equisingular stratum of the space of
quartics.
\qed
\endcorollary

We denote by~$\CM\sd$ the equisingular
stratum containing transversal determinantal
quartics. The corresponding configuration and singular homological type are
denoted by $\tS\sd$ and~$\CH\sd$, respectively; they extend the
set of singularities $\Sigma\sd:=10\bA_1$. Let
$a_1,\ldots,a_{10}\in\Sigma\sd$ be generators of the $\bA_1$
summands.

\lemma\label{Sdet}
The extension $\tS\sd\supset S\sd:=\Sigma\sd\oplus\Z h$ is
obtained from~$S\sd$ by adjoining the element
$\frac12(a_1+\ldots+a_{10}+h)$.
One has $\tS\sd\cong\bU\oplus\bE_8(2)\oplus[-4]$.
\endlemma

\proof
The requirement that $\tS$ should be an even integral lattice
implies that $\tS\sd$ is generated in $S\sd\otimes\Q$ by~$S\sd$
and several elements of the form
\roster
\item\local1
$\frac12(a_1+\ldots+a_4)$,
\item\local2
$\frac12(a_1+\ldots+a_8)$,
\item\local3
$\frac12(a_1+a_2+h)$,
\item\local4
$\frac12(a_1+\ldots+a_6+h)$,
\item\local5
$\frac12(a_1+\ldots+a_{10}+h)$.
\endroster
(up to reordering of the basis elements~$a_i$).

Case~\loccit1 is impossible as the only nontrivial finite index
extension of~$4\bA_1$ is~$\bD_4$, which contradicts
to~\iref{configuration}{conf.1}.

Consider case~\loccit2, \ie, assume that $\tS\sd$ contains
$a:=\frac12(a_1+\ldots+a_8)$. If $\tS\sd$ contained another
element~$a'$ of the same form, then, up to a further reordering,
one would have $a'=\frac12(a_1+\ldots+a_6+a_9+a_{10})$, and the
difference $a'-a$ would be as in case~\loccit1. Hence, $a$ is the
only element of \smash{$\tS\sd\bmod S\sd$} of this form, and each
surface in the stratum
has eight distinguished singular points. This contradicts
Lemma~\ref{irreducible}.

Case~\loccit3 contradicts to~\iref{configuration}{conf.2}.

Since cases~\loccit1 and~\loccit2 have been eliminated,
\smash{$\tS\sd\bmod S\sd$} may contain at most one element as
in~\loccit4 or~\loccit5. In case~\loccit4, each surface
in the stratum would have
six distinguished singular points, which would contradict
Lemma~\ref{irreducible}. Thus, either \smash{$\tS\sd=S\sd$} or
\smash{$\tS\sd\supset S\sd$}
is the index~$2$ extension generated by the
(only) element~\loccit5. Pick a point $p\in\Cp3$ and a
sufficiently generic transversal quartic $\Delta_V$, so that
$H_p\cap\Delta_V$ is nonsingular. Since $\tS\sd$ is the primitive
hull of~$S\sd$ in $L\cong H_2(\tilde\Delta_V)$, from
the Poincar\'e--Lefschetz duality it follows that
$H_1(\Delta^\circ_V\sminus H_p)=\Ext(\tS\sd/S\sd,\Z)$. Due to
Corollary~\ref{cor.pi1}, one has $[\tS\sd:S\sd]=2$,
and the first statement follows. The
isomorphism class of~$\tS\sd$ is given by a simple computation of
the discriminant group and Nikulin's uniqueness theorem
\cite{Nikulin, Theorem~1.14.2}.
\endproof

\Remark
Alternatively, case~\loccit2 in the proof of Lemma~\ref{Sdet}
can also be eliminated using
Corollary~\ref{cor.pi1}, and case~\loccit4 can be eliminated using
a refinement of this corollary stating that the group
$\pi_1(\Delta^\circ_V\sminus H_p)$ is generated by the group of
the link of any of the singular points.
\endRemark

\theorem\label{Hdet}
The configuration $\tS\sd$ given by Lemma~\ref{Sdet} extends to a
unique, up to isomorphism,
singular homological type~$\CH\sd$, which is symmetric. Thus,
one has
$\CM\sd=\CM(\CH\sd)$.
\endtheorem

\proof
The uniqueness of a primitive
embedding $\tS\sd\into L$ of the lattice~$\tS\sd$ given by
Lemma~\ref{Sdet}
follows from
\cite{Nikulin, Theorem 1.14.4}. One has
$\tS\sd^\perp\cong\bU\oplus\bE_8(2)\oplus[4]$.
Clearly, \smash{$\tS\sd^\perp$} has a vector of square~$2$, and
the reflection against the hyperplane orthogonal to such a vector
is an orientation reversing automorphism.
\endproof

\Remark
Alternatively, the fact that $\CH\sd$ is symmetric follows from
the obvious existence of real determinantal quartics.
\endRemark

\subsection{Proof of Lemma~\ref{non-coplanar}}\label{proof.non-coplanar}
We prove a stronger statement: a plane $W\subset V$ cannot contain
more than six singular points of a transversal determinantal quartic
$\Delta_V\subset V$.

Assume that seven singular points $Q_1,\ldots,Q_7$
of~$\Delta_V$ belong to a single plane $W\subset V$. Since
$\Delta_V$ is irreducible, the intersection
$\Delta_W:=\Delta_V\cap L$ is a curve, which is of degree~$4$.
Furthermore, each point $Q_1,\ldots,Q_7$ is singular
for~$\Delta_W$.

Since a reduced plane quartic has at most six singular points,
$\Delta_W$ must have multiple components.

\lemma\label{<=3}
A pencil $U\subset\Qu$ not contained entirely in~$\Delta'$
intersects~$\Delta'$ at at most three points.
\endlemma

\proof
First, assume that $U$ has a base
point singular for all quadrics. Projecting from this point, one
obtains a pencil of plane conics, singular conics corresponding to
the elements of the intersection $U\cap\Delta'$,
and the statement follows from the
fact that $\deg\Delta(2)=3$.

Now, assume that $U$ does not have a singular base point. Let
$P_1,P_2\in U\cap\Delta'$ be two distinct members of~$U$
of corank at least~$2$; they
generate~$U$.
Since $P_1$ and~$P_2$ have no common
singular points, in appropriate homogeneous coordinates one has
$P_1=\{x_0x_1=0\}$ and $P_2=\{x_2x_3=0\}$, and it is immediate
that $P_1$ and~$P_2$ are the only singular members of~$U$.
\endproof

Lemma~\ref{<=3} rules out the possibility that $\Delta_W$ contains
a double line: at most two double lines in $\Delta_W$ would
contains at most six quadrics in~$\Delta'$.

The remaining possibility is that $\Delta_W$ is a double conic. In
this case $\Delta_W$ has no linear components; in particular, no
two quadrics $P_1,P_2\in\Delta_W$ have a common singular point.
The quadric~$Q_1$ splits into two distinct
planes $\Ga_{1}$, $\Ga_{2}$.
The linear system~$W$ restricts to a pencil of conics
in~$\Ga_{1}$ containing at least six distinct singular members,
namely,
the restrictions of $Q_2,\ldots,Q_7$. (If the restrictions of two
distinct quadrics~$Q_i$, $Q_j$ coincided, $Q_i$ and~$Q_j$ would
have a common singular point.) Since $\deg\Delta(2)=3<6$,
all members of the restricted pencil are singular. Hence, the
vertex of each quadric $P\in\Delta_W\sminus\Delta'$ belongs
to~$\Ga_{1}$. The same argument shows that the vertex of~$P$ also
belongs to~$\Ga_2$, \ie, $P$ and~$Q_1$ have a common singular
point. This contradiction concludes the proof of
Lemma~\ref{non-coplanar}.
\qed

\Remark
Using the surjectivity of the period map and the Riemann--Roch
theorem for $K3$-surfaces, one can easily show that
the stratum $\CM\sd$ does contain a quartic with all ten singular
points coplanar (lying in a curve of degree two).
In particular, determinantal quartics form a
proper subset of~$\CM\sd$.
\endRemark

\section{Real determinantal quartics\label{S.R}}

In this section, we discuss the topology of a determinantal
quartic with nonempty spectrahedral domain and
prove
Theorem~\ref{th.main}.

\subsection{Geometric real structures}
We always consider the space $\Qu(n)$ with its
\emph{geometric real structure}, \ie, the one induced by
the complex conjugation $\conj\:\Cp{n}\to\Cp{n}$.
All real quadratic forms
%form
constitute
a linear space $\R^{N(n)+1}$, and
one has a double covering $S^{N(n)}\to\Qu(n)_\R$, where
$S^{N(n)}\subset\R^{N(n)+1}$ is the unit sphere.
%In the sequel, we fix
%a basis in $\R^{N(n)+1}$ and identify $S^{N(n)}$ with the unit
%sphere.

We reserve the notation~$\spbar$ for the lift from $\Qu(n)_\R$ to
$S^{N(n)}$; in particular, one has real discriminant hypersurfaces
$\bar\Delta(n)\subset S^{N(n)}$ and $\bar\Delta\subset S^9$.

Recall that a real quadratic form~$\bar q$ has a well defined \emph{index}
$\ind\bar q$
(the negative inertia index of~$\bar q$); one has
$0\le\ind\bar q\le n+1$.

A \emph{real determinantal quartic} is a real quartic
$X\subset\Cp3$ equivariantly isomorphic to a quartic
$\Delta_V\subset V$, where $V\subset\Qu$ is a $3$-subspace real
with respect to the geometric real structure. Given such a
quartic~$X$,
the
\emph{spectrahedral region} of~$X$ is the (only) connected
component of the complement $\Rp3\sminus X_\R$ constituted by the
quadrics represented by quadratic forms of index~$0$
(equivalently, those of maximal index~$4$).

\lemma\label{x4}
Let $X\subset\Cp3$ be a real determinantal quartic.
Then any real line
meeting the spectrahedral region of~$X$
intersects~$X$ at four real points \rom(counted with
multiplicities\rom). In other words, all intersection points are
real.
\endlemma

\proof
Identify $(\Cp3,X)$ with a real pair $(V,\Delta_V)$ and let~$W$
be the image of the line in question. Consider the lift
$\bar W\subset S^9$. The index function $\ind\:\bar W\to\Z$ is
locally constant on $\bar W\sminus\bar\Delta$, and its
increment~$\delta_p$ at an intersection point
$p\in\bar W\cap\bar\Delta$ of multiplicity~$m_p$ is subject to the
conditions $\mathopen|\delta_p\mathclose|\le m_p$,
$\delta_p=m_p\bmod2$. Any point $Q$ in the spectrahedral region
of~$\Delta_V$ lifts to a pair~$\bar q$, $-\bar q$ of
quadratic forms, so that one has $\ind\bar q=4$,
$\ind(-\bar q)=0$. Since the two indices differ by~$4$, the
intersection $\bar W\cap\bar\Delta$ must consist of at least eight
points (counted with multiplicities).
\endproof

\corollary\label{2.spheres}
Let $X\subset\Cp3$ be a real determinantal quartic with nonempty
spectrahedral region. Then~$X$ is a non-positive degeneration of a
nonsingular real quartic $X'\subset\Cp3$ with the real part
$X'_\R$ constituted by a nested pair of spheres. \qed
\endcorollary

\subsection{Preliminary computation}\label{s.prelim}
Consider the
%nonsingular
real homological type
$(L, h, c)$
corresponding to
%a nonsingular real quartic $X\subset\Rp3$
nonsingular real quartics
with two nested spheres; such quartics are amphichiral.
According
to~\cite{Nikulin},
one has
$$
L^c_+\cong\bE_8(2)\oplus2\bA_1\oplus\bU,\qquad
L^c_-\cong\bE_8(2)\oplus2\bA_1(-1).
\eqtag\label{eq.L}
$$
Fix standard bases $e_1,\ldots,e_8$, $v_1$, $v_2$, and
$u_1$, $u_2$ for $\bE_8(2)$, $2\bA_1$, and~$\bU$, respectively,
and let $e'_1,\ldots,e'_8$, $v'_1$, $v'_2$ be the `matching'
standard basis for~$L^c_-$, so that the sum $r+r'$ of two basis
vectors $r\in L^c_+$, $r'\in L^c_-$ of the same name is divisible
by~$2$ in~$L$.
In view of~\cite{Nikulin},
$$
h=v'_1+v'_2,\qquad
h=v_1+v_2\bmod2L.
\eqtag\label{eq.h}
$$

Let $\fS\subset\fP\subset\bar\fP$ be a
fundamental tower
of $(L, h, c)$.
The
polyhedron $\fS$ is finite; its Coxeter scheme, computed
in~\cite{Student}, is shown in Figure~\ref{fig.poly}, and $\fS$
can be chosen to be bounded by
the walls orthogonal to the vectors indicated in the
figure. (In the figure, walls of type~\iref{def.poly}{s2}
and~\ditto{s4} are represented by $\circ$
and~$\bullet$, respectively, and the only wall of
type~\iref{def.poly}{s4h} is represented by a circled bullet.
Whenever
the hyperplanes supporting
two walls
intersect at an angle $\pi/n$, $n\ge2$, the corresponding
vertices  are connected by $(n-2)$ edges.
Note that, in the case under consideration, any two walls
do intersect.)

Clearly, $\fP$ is paved by the (infinitely many) copies of~$\fS$
obtained by iterated reflections against the walls of
type~\iref{def.poly}{s4} (vertices $e_1,\ldots,e_8$, $e_{12}$,
and~$e_{13}$
in the figure), and
$\bar\fP$ is the union of~$\fP$ and its image under the reflection
against the only wall~$e_9$ of type~\iref{def.poly}{s4h};
%(vertex $e_9$);
see Definition~\ref{def.poly} and Remark~\ref{rem.fS}.

\midinsert
%$$
\smallskip
\centerline{%
\DDlinelength=6
\DDbox
@(\bullet)@*<|><1|>@---
 @(\bullet)@*<|><2|>@---
 @(\bullet)@*<|><3|>@---
 @(\bullet)@*<|><5|>@---
 @(\bullet)@*<|><6|>@---
 @(\bullet)@*<|><7|>\cr
@|---@.@|---@.@.@|---\cr
@(\bullet)@*<13|><|>@.@.@(\bullet)@*<4|><|>@.@.@.
 @(\bullet)@*<|8><|>\cr
@|===@.@.@.@.@|---\cr
@(\circ)@*<|><|10>@===
 @(\Circ\bullet)@*<|><|9>\expand3@===
 @(\circ)@*<|><|11>@---
 @(\circ)@*<|><|0>@===
 @(\bullet)@*<|><|12>
\endDD}
%$$
$$
\gather
\alignedat2
&e_0&&=u_1-u_2,\\
&e_9&&=v_1-v_2,\\
&e_{10}&&=v_2,
\endalignedat\qquad
\alignedat2
&e_{11}&&=u_2-v_1,\\
&e_{12}&&=2u_2+e_8^*,\\
&e_{13}&&=2(u_1+u_2)-v_1-v_2+e_1^*,
\endalignedat\\
\noalign{\medskip}
\aligned
e_1^*&=-4e_1-7e_2-10e_3-5e_4-8e_5-6e_6-4e_7-2e_8,\\
e_8^*&=-2e_1-4e_2-6e_3-3e_4-5e_5-4e_6-3e_7-2e_8.
\endaligned
\endgather
$$
\par\removelastskip
%\smallskip
\figure\label{fig.poly}
The fundamental polyhedron~$\fS$
\endfigure
\endinsert

Let $X\subset\Cp3$ be a properly marked nonsingular quartic of
type $(L,h,c)$.

\lemma\label{spheres}
The classes realized in~$L^c_+$
by the inner and outer spheres
of~$X_\R$ are $\sin=e_{11}$ and $\sout=e_{11}+e_9$, respectively.
\endlemma

\proof
Let $G$ be the graph obtained from the Coxeter scheme of~$\bar\fP$
by removing all but simple edges.
According to~\cite{DIK, Theorem 16.1.1},
any vertex of~$G$ of valency~$>2$
is a class realized by a spherical component of~$X_\R$.
Clearly, $e_{11}$ and $e_{11}+e_9$ (obtained
from~$e_{11}$ by reflection) are two such vertices.
%Recall, see~\cite{DIK},
%that any vertex of valency~$>2$
%of the Coxeter scheme
%of~$\bar\fP$ is a class realized by a spherical component
%of~$X_\R$, and $e_{11}$ and $e_{11}+e_9$ (obtained
%from~$e_{11}$ by reflection) are clearly two such vertices.
The
outer sphere is definitely not contractible. Hence, the
vector~$e_{11}$, which defines a wall of~$\fP$ and thus can serve
as a vanishing cycle, see Theorem~\ref{th.if},
represents the inner sphere.
\endproof

Denote by $L_+^0\subset L^c_+$ the sublattice spanned by
$e_1,\ldots,e_8$, $e_{12}$, and~$e_{13}$.

\lemma\label{vanishing}
Each real vanishing cycle $r\in L^c_+$
of a non-positive nodal
degeneration of~$X$ is of one of the following three
forms\rom:
\roster
\item\local{v.inner}
$e_{11}=\sin$ \rom(the inner sphere shrinks to a point\rom)\rom;
\item\local{v.both}
$e_0+d$, $d\in L_+^0$ \rom(a common point of the two
spheres\rom)\rom;
\item\local{v.outer}
$e_{10}+d$, $d\in L_+^0$ \rom(a node in the outer sphere\rom).
\endroster
For each pair $s'$, $s''$ of conjugate vanishing cycles, the
invariant decomposable vector $s'-s''$ belongs to~$L_+^0$.
\endlemma

\proof
Each real vanishing cycle is a wall of~$\fP$ of
type~\iref{def.poly}{s2}, see Theorem~\ref{th.only.if}.
From the description of~$\fP$ in terms
of~$\fS$ and Figure~\ref{fig.poly} it follows that any such wall
either is~$e_{11}$ or is obtained from~$e_0$ or~$e_{10}$ by
iterated reflections against the walls
$e_1,\ldots,e_8$, $e_{12}$, and~$e_{13}$ (and their consecutive images). The
geometry of the corresponding degeneration is easily seen from
comparing the vanishing cycle~$r$ against the classes of the spheres:
one has $r=\sin$ in case~\loccit{v.inner},
$r\cdot\sin=r\cdot\sout=1$ in case~\loccit{v.both}, and
$r\cdot\sin=0$, $r\cdot\sout=2$ in case~\loccit{v.outer}.

For a pair $s'$, $s''$ of conjugate vanishing cycles, the vector
$s'-s''$ is either one of the type~\iref{def.poly}{s4} walls
of~$\fS$ or one of their consecutive images
under reflections, see Theorem~\ref{th.only.if}.
\endproof

\lemma\label{no.5}
The lattice~$L^c_-$
%cannot
does not
contain a quintuple~$t_i$, $i=1,\ldots,5$,
of pairwise
orthogonal vectors of square~$(-4)$ such that
$t_1+\ldots+t_5=h\bmod2L^c_-$.
\endlemma

\proof
In view of~\eqref{eq.L} and~\eqref{eq.h},
the vector~$h$ is
characteristic in the sense that $a^2+a\cdot h=0\bmod4$ for any
$a\in L^c_-$. Hence, $(h+2a)^2=4\bmod16$ for any $a\in L^c_-$.
On the other hand, $(t_1+\ldots+t_5)^2=-20=-4\bmod16$.
\endproof

\subsection{Proof of Theorem~\ref{th.main}}\label{proof.main}
We keep the notation of Subsection~\ref{s.prelim}.

The assumption that the spectrahedral region~$R$ of~$X$
is nonempty rules out
real vanishing cycles of type~\iref{vanishing}{v.inner}. Assume
that $X$ has $m$ vanishing cycles of type~\iref{vanishing}{v.both}
and $n$ vanishing cycles of type~\ditto{v.outer} (and $10-m-n$
imaginary vanishing cycles split into conjugate pairs). Using
Lemma~\ref{vanishing} and the description of the vectors involved
given in Figure~\ref{fig.poly}, one can easily see that the
parities of the coefficients of~$v_1$ and~$v_2$ in the sum of all
ten vanishing cycles differ by $n\bmod2$. Due to~\eqref{eq.h} and
Lemma~\ref{Sdet}, $n$ is even, and so is~$m$.

If $m=n=0$, then~$X$ has five pairs of complex conjugate vanishing
cycles $s'_j$, $s''_j=-c(s'_j)$, $j=1,\ldots,5$,
and the skew-invariant vectors
$s'_j+s''_j\in L^c_-$
form a quintuple contradicting Lemma~\ref{no.5}.

For the construction, relabel the nine vertices of
type~\iref{def.poly}{s4} in the edges of the Coxeter scheme
consecutively, \ie, let $e_{13}=w_1$, $e_i=w_{i+1}$ for $i=1,2,3$,
$e_i=w_i$ for $i=5,\ldots,8$, and $e_{12}=w_9$. Pick
a pair $m$, $n$ of even integers as
in the statement, denote $p=5-\tfrac12(m+n)$,
and consider the following vectors:
$$
\alignedat2
&r'_i=e_0+w_9+\ldots+w_{11-i},&\qquad
 &\text{$i=1,\ldots,m$ (if $m>0$)},\\
&r''_j=e_{10}+w_1+\ldots+w_{j-1},&\qquad
 &\text{$j=1,\ldots,n$ (if $n>0$)},\\
&t_k=w_{n+2k-1},&\qquad
 &\text{$k=1,\ldots,p$ (if $p>0$)}.
\endalignedat
$$
It is straightforward to check that:
\roster
\item
each~$r'_i$ is obtained by a sequence of reflections from
the vertex~$e_0$,
\ie, is as in Lemma~\iref{vanishing}{v.both};
\item
each~$r''_j$ is obtained by a sequence of reflections from
the vertex~$e_{10}$,
\ie, is as in Lemma~\iref{vanishing}{v.outer};
\item
each~$t_k$ is a wall of~$\fS$ of type~\iref{def.poly}{s4};
\item
all vectors are pairwise orthogonal;
\item
all vectors are linearly independent in $L^c_+/2L^c_+$;
\item
the sum of all ten vectors equals $h\bmod2L$.
\endroster
Furthermore,
assuming that $p\le4$, one can easily find pairwise orthogonal
vectors $t'_1,\ldots,t'_p\in L^c_-$ such that $t^{\prime2}_k=-4$,
$t'_k\cdot h=0$, and $t'_k=t_k\bmod2L$, $k=1,\ldots,p$.
Indeed,
consider the vectors $w_1'=v'_1-v'_2+e_1^{\prime*}$,
$w_3'=e_2'$, $w_5'=e_5'$, $w_7'=e_7'$, $w_9'=e_8^{\prime*}$.
%These vectors
They
are orthogonal to~$h$ and have square~$(-4)$, and
any sequence of up to four consecutive vectors is orthogonal. (In
fact, all five vectors are pairwise orthogonal except that
$w'_1\cdot w'_9\ne0$.)
Now, one can take
for $t_k'$ the `matching' vectors $w'_*$.

Finally, the set~$\Gs$ constituted by the ten vectors
$$
r'_1,\ldots,r'_m,\ r''_1,\ldots,r''_n,\
 \tfrac12(t'_1\pm t_1),\ldots,\tfrac12(t'_p\pm t_p)
$$
is an admissible system of cycles, see Definition~\ref{def.cycles};
due to Theorem~\ref{th.if}, it
can serve as the set of vanishing cycles of a
non-positive nodal degeneration of~$X$.
On the other hand, the set~$\Gs$
satisfies
Lemma~\ref{Sdet}; hence,
according to Corollary~\ref{whole.stratum} and
Theorem~\ref{Hdet},
a generic degeneration of~$X$ contracting these
vanishing cycles is a transversal determinantal
quartic.
\qed

\refstyle{C}

\break

\enddocument